\documentclass{elsart}
\usepackage{amscd,amssymb,hyperref}
\usepackage{hyperref}
\usepackage{subfigure}
\usepackage[graph,frame,poly,arc,matrix]{xy}  

\newtheorem{example}{Example}[section]
\newtheorem{question}{Question}[section]

\newtheorem{mthm}{Theorem} 

\newcommand{\abs}[1]{\left|#1\right|}
\renewcommand{\b}[1]{\mathbf{#1}}

\newcommand{\A}{\mathcal{A}}
\newcommand{\B}{\mathcal{B}}

\newcommand{\DD}{\mathcal{D}}
\newcommand{\FF}{\mathcal{F}}

\newcommand{\NN}{\mathcal{N}}
\newcommand{\RR}{\mathcal{R}}

\newcommand{\ZZ}{\mathcal{Z}}
\newcommand{\W}{\mathcal{W}}
\newcommand{\cV}{{\mathcal{V}}}
\newcommand{\wA}{\widehat{\mathcal{A}}}

\newcommand{\C}{\mathbb{C}}
\newcommand{\F}{\mathbb{F}}
\newcommand{\Z}{\mathbb{Z}}
\newcommand{\R}{\mathbb{R}}

\renewcommand{\a}{\alpha}

\newcommand{\Hom}{\operatorname{Hom}} 
 
\newcommand{\rank}{\operatorname{rank}}

\newcommand{\id}{\operatorname{id}}
\newcommand{\Aut}{\operatorname{Aut}}

\newcommand{\GL}{\operatorname{GL}}
\newcommand{\ab}{\operatorname{ab}}

\newcommand{\gr}{\operatorname{gr}}

\def\Re{\operatorname{Re}}

\newcommand{\ii}{i}

\begin{document}

\begin{frontmatter}

\title
{Translated tori in the characteristic varieties of complex 
hyperplane arrangements}

\author{Alexander~I.~Suciu\thanksref{email}}
\address{Department of Mathematics,
Northeastern University,
Boston, MA 02115}

\thanks[email]{E-mail: \href{mailto:alexsuciu@neu.edu}{alexsuciu@neu.edu}. 
Research supported by Northeastern University (through an RSDF
grant),  and by the Volkswagen-Stiftung (RiP-program at Oberwolfach).}




\begin{abstract}
We give examples of complex hyperplane arrangements $\A$  
for which the top characteristic variety, $V_1(\A)$,  
contains positive-dimensional irreducible components 
that do not pass through the origin of the algebraic 
torus $(\C^*)^{|\A|}$.  These examples answer several 
questions of Libgober and Yuzvinsky. As an application, 
we exhibit a pair of arrangements for which 
the resonance varieties of the Orlik-Solomon algebra 
are (abstractly) isomorphic, yet whose characteristic 
varieties are not isomorphic.  The difference comes from 
translated components, which are not detected by the 
tangent cone at the origin.
\end{abstract}

\end{frontmatter}

\section{Introduction}

The characteristic varieties of a space $X$ are the jumping loci 
of the cohomology of $X$ with coefficients in rank $1$ local systems:
$
V_d(X)=\{ \b{t}\in (\C^*)^{b_1(X)} \mid \dim_{\C} H^1(X, \C_{\b{t}})\ge d\}.
$
If $X$ is the complement of a normal-crossing divisor in a compact 
K\"{a}hler manifold with vanishing first homology, then $V_{d}(X)$ is 
a finite union of torsion-translated subtori of the character torus, 
see \cite{Ar}.  If $\A$ is an arrangement of hyperplanes in $\C^{\ell}$, 
with complement $X=X(\A)$,  
then the irreducible components of $V_{d}(\A):=V_{d}(X)$ which  
contain the origin can be determined combinatorially, from the intersection 
lattice of $\A$.  This follows from the fact that the tangent cone to
$V_{d}(\A)$ at $\b{1}$ coincides with the resonance variety, $\RR_{d}(\A)$, 
of the Orlik-Solomon algebra (see \cite{CScv} for a proof, and 
\cite{Li3}, \cite{CO}, \cite{Li4} for other proofs and generalizations).  
The variety $\RR_{d}(\A)$ in turn admits explicit combinatorial descriptions, 
see \cite{Fa}, \cite{LY1}. 

It was first noted in \cite{CScv} that there exists a hyperplane arrangement 
for which $V_{2}$ contains translated tori.  These translated tori 
are isolated torsion points in $V_{2}$, lying at the intersection of several 
components of $V_{1}$ which do pass through the origin 
(see Example~4.4 in \cite{CScv} and Example~\ref{ex:diamond} below).   
Thus, the question arose whether the characteristic varieties of a 
complex hyperplane arrangement may have {\it positive-dimensional} 
translated components, see \cite{LY2}, Problem~5.1.  In this note, 
we answer that question, as follows.  

\begin{mthm}
\label{thm:transtorus}
There exist arrangements of complex hyperplanes for 
which the top characteristic variety, $V_{1}$, contains 
positive-dimensional irreducible components which do 
not pass through the origin. 
\end{mthm}

\begin{figure}
\setlength{\unitlength}{0.5cm}
\begin{picture}(16,10)(-5,0)
\put(1,8){\line(1,0){14}}
\put(14,1){\line(0,1){8}}
\put(2,8){\line(3,-1){13}}
\put(2,8){\line(-3,1){1.15}}
\put(2,8){\line(4,-1){13}}
\put(2,8){\line(-4,1){1.1}}
\put(8,8){\line(1,-1){7}}
\put(8,8){\line(-1,1){0.8}}
\put(8,8){\line(3,-2){7}}
\put(8,8){\line(-3,2){1.1}}
\put(10,8){\line(2,-3){4.65}}
\put(10,8){\line(-2,3){0.5}}
\put(14,8){\line(-1,-1){4.5}}
\put(14,8){\line(1,1){0.9}}
\put(15.7,9){\makebox(0,0){$-1$}}
\put(15.5,8){\makebox(0,0){$t$}}
\put(15.9,4.9){\makebox(0,0){$t^{-2}$}}
\put(15.5,3.7){\makebox(0,0){$t$}}
\put(15.7,3){\makebox(0,0){$-1$}}
\put(16,1.2){\makebox(0,0){$-t^{-1}$}}
\put(14.8,0.5){\makebox(0,0){$t^{2}$}}
\put(13.0,0.8){\makebox(0,0){$-t^{-1}$}}
\end{picture}
\caption{\textsf{Generic section of deleted ${\rm B}_3$ 
arrangement and translated torus in $V_1$}}
\label{fig:deletedb3}
\end{figure}

The simplest such arrangement is the ``deleted ${\rm B}_3$" 
arrangement, $\DD$, discussed in Example~\ref{ex:deletedb3}.
It is an arrangement of $8$ planes in $\C^3$, 
with defining forms $x-z$, $y-z$, $x$, $y$, $x-y+z$, $z$, 
$x-y-z$, $x-y$.   The arrangement $\DD$ is fiber-type, with 
exponents $\{1,3,4\}$. The variety $V_1(\DD)$ has a component 
parametrized by $\{(t,-t^{-1},-t^{-1},t,t^2,-1,t^{-2},-1) \mid  
t\in \C^*\}$.  This is a $1$-dimensional torus, translated by 
a second root of $\b{1}$.   Figure~\ref{fig:deletedb3} depicts 
the real part of a generic $2$-dimensional section of $\DD$ 
(obtained by setting $z=2x+3y+1$), 
together with the local system corresponding to the point $t\in \C^*$. 

The deleted ${\rm B}_3$ arrangement can also be used to answer 
Conjecture 4.4 from \cite{LY1}, and Problems 5.2 and 5.3 from \cite{LY2}. 

As noted in \cite{Li3}, \cite{LY1} (see also \cite{Fa}), 
all the components of $V_1$ passing through the origin 
must have dimension at least $2$.  On the other hand, 
all the positive-dimensional translated components 
that we find correspond to $\DD$ sub-arrange\-ments, 
and so have dimension $1$.   We do not know
whether translated components can have dimension greater than $1$, 
but we exhibit an arrangement where $V_{1}$ has $0$-dimensional 
components (Example~\ref{ex:grunbaum}).  

The translated components in the characteristic varieties of an 
arrangement $\A$ are not detected by the tangent cone at the origin, 
and thus contain information which is not available from the 
Orlik-Solomon algebra of $\A$, at least not directly.  
We illustrate this phenomenon in Example~\ref{ex:ziegler}, where 
we find a pair of arrangements for which the resonance varieties 
are (abstractly) isomorphic,  but the characteristic varieties 
have a different number of components.  These two arrangements 
have non-isomorphic lattices, though.  Thus, it is still an open 
question whether the translated components of $V_d(\A)$ are 
combinatorially determined.

One of the main motivations for the study of characteristic varieties 
of a space $X$ is the very precise information they give about the 
homology of finite abelian covers of $X$, see \cite{Li2}, \cite{Sa}.   
From that point of view, the existence of translated components in 
$V_{1}$ has immediate repercussions on the Betti numbers of some 
finite covers of $X$ (those corresponding to torsion characters belonging 
to that component).  But it also affects the torsion coefficients of 
some abelian covers of $X$, and the number of certain metabelian covers 
of $X$.  These aspects are pursued in joint work with D.~Matei, \cite{MS3}.  
The starting point of that paper was the discovery of $2$-torsion in 
the homology of certain $3$-fold covers of $X(\DD)$.  We were led to the 
translated component in $V_1(\DD)$ by an effort to explain that 
unexpected torsion.

\section{Characteristic varieties and hyperplane arrangements}
\label{sec:cvar}

We start by reviewing methods for computing the fundamental group, 
the characteristic varieties, and the resonance varieties of a 
complex hyperplane arrangement.

\subsection{Characteristic varieties}
\label{subsec:cv}
Let $X$ be a space having the homotopy type of a connected, 
finite CW-complex.   For simplicity, we will assume throughout 
that $H=H_1(X,\Z)$ is torsion free.  Set $n=b_1(X)$, and fix 
a basis $\{t_1,\dots,t_n\}$ for $H\cong\Z^n$.  Let $G=\pi_1(X)$ be 
the fundamental group, and $\ab:G\to H$ the 
abelianization homomorphism.

Let $\C^*$ be the multiplicative group of units in $\C$, and let
$\Hom(G,\C^*)$  be the group of characters of $G$.   
Notice that $\Hom(G,\C^*)$ is isomorphic to the affine algebraic group
$\Hom(H,\C^*)\cong (\C^*)^n$, with coordinate ring 
$\C H\cong\C[t_1^{\pm{1}},\dots, t_n^{\pm{1}}]$.   
For each integer $d\ge 0$, set 
\[
V_d(X)=\{\b{t}=(t_1,\dots ,t_n)\in (\C^{*})^n \mid 
\dim_{\C} H^1(G,\C_{\mathbf{t}}) \ge d\},
\]
where $\C_{\b{t}}$ is the $G$-module $\C$ 
given by the representation 
$\begin{xy}\xymatrix@1{G \ar[r]^{\ab} &\Z^n \ar[r]^{\b{t}} &\C^*}\end{xy}$. 

Then $V_d(X)$ is an algebraic subvariety of the complex $n$-torus, called the 
{\em $d$-th characteristic variety} of $X$.  
The characteristic varieties form a descending tower,  
$(\C^{*})^n=V_0\supseteq V_{1}\supseteq \cdots \supseteq V_{n-1}\supseteq
V_{n}$, which depends only on the
isomorphism type of $G=\pi_1(X)$, up to a monomial 
change of basis in $(\C^*)^n$, see \cite{MS1}.   

As shown in~\cite{Hr}, the characteristic varieties of $X$ may be interpreted 
as the determinantal varieties of the Alexander matrix of the group $G=\pi_1(X)$. 
Given a presentation $G=\langle x_1,\dots , x_m\mid r_1,\dots , r_s\rangle$, the 
Alexander matrix is the $s \times m$ matrix 
$A=\left(\frac{\partial r_i}{\partial x_j}
\right)^{\ab}$, with entries in $\C[t_1^{\pm{1}},\dots, t_n^{\pm{1}}]$, 
obtained by abelianizing the Jacobian of Fox derivatives of the relations.  Let
$A(\b{t})$ be the evaluation of $A$ at $\b{t}\in (\C^{*})^n$.   
For $0\le d< n$, we have
\[
V_d(X)=\{\b{t}\in (\C^*)^n \mid \rank A(\b{t}) < m-d\}. 
\] 

Remarkably, the existence of certain analytic or geometric structures on a space  
puts strong qualitative restrictions on the nature of its characteristic varieties. 
There are several results along these lines,  due to Green, Lazarfeld, Simpson, and
Arapura.  The result we need is  the following:

\begin{mthm}[Arapura~\cite{Ar}] 
\label{thm:arapura}  
Let $X$ be the complement of a normal-crossing divisor in a compact 
K\"{a}hler manifold with vanishing first homology.  Then each characteristic 
variety $V_{d}(X)$ is a finite union of torsion-translated subtori 
of the algebraic torus $(\C^*)^{b_1(X)}$.
\end{mthm}

\subsection{Fundamental groups of arrangements}
\label{subsec:pi1arr}
Let $\A =\{H_{1},\dots ,H_{n}\}$ be an arrangement of (affine) 
hyperplanes in $\C^{\ell}$, $\ell\ge 2$, with complement 
$X(\A)=\C^{\ell}\setminus \bigcup_{i=1}^n H_{i}$.  
We review the procedure for finding the braid monodromy 
presentation of the fundamental group of the complement, 
$G(\A)=\pi_1(X(\A))$. 
This presentation is equivalent to the Randell-Arvola 
presentation (see \cite{CSbm}), and the $2$-complex 
modelled on it is homotopy-equivalent to $X=X(\A)$ (see \cite{Li1}).   
Since we are only interested in $G=G(\A)$, the well-known 
Lefschetz-type theorem of Hamm and L\^{e} allows us to assume 
that $\ell=2$, by replacing $\A$ with a generic $2$-dimensional 
slice, if necessary.  

Let $v_{1},\dots, v_{s}$ be the intersection points of the 
lines of $\A$. The combinatorics of the
arrangement is encoded in its intersection poset,
$L(\A)=\{L_1(\A),L_2(\A)\}$, where $L_1=\b{n}:=\{1,\dots , n\}$ 
and $L_2=\{I_{1},\dots,I_{s}\}$, with 
$I_{k}=\{i \in \b{n} \mid H_{i}\cap v_{k}\ne \emptyset\}$. 
Choosing a generic linear projection $p:\C^2\to \C$, and a 
basepoint $y_0\in \C$ such that 
$\Re(y_0)>\Re(p(v_1))>\cdots>\Re(p(v_s))$
gives orderings of the lines and vertices, which we may 
assume coincide with the orderings specified above.  
Choosing also a path in $\C$, starting at
$y_0$,  and passing successively through $p(v_1),\dots,p(v_s)$ 
gives a ``braided wiring diagram,"   
$\W(\A)=\{I_{1},\beta_{1},I_{2},\dots,\beta_{s-1},I_{s}\}$, 
where $\beta_{k}$ are certain braids in the Artin braid group $B_n$. 

Let $\{A_{i,j}\}_{1\le i<j\le n}$ be the usual generating set 
for the pure braid group $P_n$, as specified in \cite{Bi}.  
More generally, for $I\subset \b{n}$,  let $A_{I}\in P_n$ 
be the full twist on the strands indexed by $I$.  The braid 
monodromy presentation of $G=\pi_1(X)$ is given by:  
\begin{equation} 
\label{eq:bmpres}
G=\langle x_1,\dots ,x_n \mid \alpha_{k}(x_i)=x_i  
\ {\rm for }\ i \in I_{k}\setminus \{\max I_{k}\}
\ {\rm and }\ k\in \b{s}\rangle,   
\end{equation}
where each $\alpha_{k}$ is a pure braid of the form 
$A_{I_{k}}^{\delta_{k}}=\delta_{k}^{-1}A_{I_{k}}\delta_{k}$, 
acting on $\F_n=\langle x_1,\dots , x_n\rangle$ 
by the Artin representation $P_n\hookrightarrow \Aut(F_n)$.  The conjugating 
braids $\delta_{k}$ may be obtained from $\W$, as follows. 

In the case where $\A$ is the complexification of a real arrangement, 
$\W$ may be realized as a (planar) wiring diagram (with all $\beta_{k}=1$), 
in the obvious way.  Each vertex set $I_{k}\in \W$  gives rise to a partition
$\b{n} = I'_{k}\cup I_{k} \cup I''_{k}$ into lower, middle, and upper wires.  
Set $J_{k}=\{i \in I''_{k} \mid \min I_{k}<i<\max I_{k}\}$.  
Then $\delta_{k}$ is the subword of 
$A_{\b{n}}=\prod_{i=2}^n\prod_{j=1}^{i-1} A_{j,i}$, 
given by $\delta_{k} = \prod_{i\in I_k} \prod_{j\in J_k} A_{j,i}$, see
\cite{CF} (and also \cite{CSbm}).   In the general case, the braids 
$\beta_1,\dots,\beta_{k-1}$ must also be taken into account, 
see \cite{CSbm} for details. 

\subsection{Characteristic and resonance varieties of arrangements}
\label{subsec:cvarr}
For an arrangement $\A$, with complement $X=X(\A)$, let $V_{d}(\A):=V_{d}(X)$.  
In equations, $V_d(\A)=\{\b{t}\in (\C^*)^n\mid \rank A(\b{t}) < n-d\}$,  
where $A$ is the Alexander matrix corresponding to the presentation 
(\ref{eq:bmpres}) of $G=\pi_1(X)$.  By Arapura's Theorem, 
$V_{d}(\A)$ is a finite union of torsion-translated 
tori in $(\C^*)^n$.  Denote by $\check{V}_{d}(\A)$ the union of those tori 
that pass through $\b{1}$, and by $\cV_{d}(\A)$ the tangent cone of 
$\check{V}_{d}(\A)$ at $\b{1}$.  
Clearly, $\cV_{d}(\A)$ is a central arrangement of subspaces in $\C^n$. 
The exponential map,  
$\exp:{\rm T}_{\b{1}}((\C^*)^n)=\C^n\to (\C^*)^n$, 
$\lambda_i \mapsto e^{2\pi\ii \lambda_i}=t_i$, takes   
each subspace in $\cV_{d}(\A)$ to the corresponding
subtorus in $\check{V}_{d}(\A)$. In equations, 
$\cV_d(\A)=\{\lambda\in \C^n\mid \rank A^{(1)}(\lambda) < n-d\}$,  
where $A^{(1)}$ is the linearized Alexander matrix of $G$, 
see \cite{CScv} (and also \cite{MS2}).   The variety
$\cV_{d}(\A)$  (and thus, $\check{V}_{d}(\A)$, too) admits a
completely combinatorial  description, as follows.

The {\it $d^{\rm{th}}$ resonance variety} of a space $X$ is the set 
$\RR_d(X)$ of cohomology classes $\lambda \in H^1(X,\C)$ for which 
there is a subspace $W \subset H^1(X,\C)$, of dimension $d+1$, such 
that $\lambda\cup W= 0$ (see \cite{MS2}).   In other words, 
\[\RR_d(X)=\{\lambda \mid \dim H^1(H^*(X,\C),\lambda)\ge d\}.\]
The resonance varieties of an arrangement, $\RR_d(\A):=\RR_d(X(\A))$, 
were introduced and studied in \cite{Fa}.   
It turns out that $\cV_{d}(\A)=\RR_d(\A)$, 
see \cite{CScv}, \cite{Li3} for two different proofs, 
and \cite{CO}, \cite{Li4} for recent generalizations. 

As seen above, the top resonance variety is the union of a subspace 
arrangement:  $\RR_1(\A)=C_1\cup\cdots\cup C_r$.  It is also known that 
$\dim C_i\ge 2$, $C_i\cap C_j=\{\b{0}\}$ for $i\ne j$, and
$\RR_{d}(\A)=\{\b{0}\}\cup\bigcup_{\dim C_i\ge d+1} C_i$, see \cite{LY1}. 
For each $I\in L_2(\A)$ with $\abs{I}\ge 3$, there is a {\em local}
component, $C_{I}=\{\lambda \mid \sum_{i}\lambda_i=0\ {\rm and}\ 
\lambda_i=0\ {\rm for}\ i\notin I\}$.  Note that $\dim C_I=\abs{I}-1$, 
and thus $C_I\subset \RR_{\abs{I}-2}(\A)$.  

The non-local components also admit a description purely in 
terms of $L(\A)$, see \cite{Fa}, \cite{LY1}.  A partition 
$\mathsf{P}=(\mathsf{p}_1\, | \cdots |\,\mathsf{p}_q)$
of $\bf{n}$ is called {\em neighborly} if, for all $I\in L_2(\A)$,
the following holds: 
$
\abs{\mathsf{p}_j\cap I}\ge \abs{I}-1 \Longrightarrow
I\subset \mathsf{p}_j.
$
To a neighborly partition $\mathsf{P}$, there corresponds a subspace
\[
C_{\mathsf{P}}=
\{\lambda \mid \sum_i\lambda_i=0\} \cap 
\bigcap_{I}\,\{\lambda \mid \sum_{i\in I}\lambda_i = 0\},
\]
where $I$ ranges over all vertex sets not contained in a single block of
$\mathsf{P}$.  Results of \cite{LY1} imply that, if $\dim C_{\mathsf{P}}
\ge 2$, then $C_{\mathsf{P}}$ is a component of $\RR_1(\A)$. All the
components of $\RR_1(\A)$ arise in this fashion from neighborly partitions of
sub-arrangements of $\A$.  

This completes the combinatorial description of
$\cV_{d}(\A)=\RR_d(\A)$, and thus, that of $\check{V}_d(\A)$.

\subsection{Decones and linearly fibered extensions}
\label{subsec:decext}
We conclude this section with two constructions which simplify in many 
instances the computation of the characteristic varieties of an arrangement.  

The first construction associates to a central arrangement  
$\A=\{H_1,\dots, H_n\}$ in $\C^{\ell}$, an affine arrangement, 
$\A^*$, of $n-1$ hyperplanes in $\C^{\ell-1}$, called a 
{\em decone} of $\A$.  Let  $Q$ be a defining polynomial for $\A$.  
Choose coordinates
$(z_1,\dots,z_{\ell})$  in $\C^{\ell}$ so that $H_n = \ker(z_{\ell})$.  
Then, $Q^* = Q(z_1,\dots,z_{\ell-1},1)$ is a defining 
polynomial for $\A^*$, and $X(\A) \cong X(\A^*) \times\C^*$, 
see \cite{OT}.   It follows that:
\[
V_{d}(\A) = \{\b{t}\in (\C^*)^n \mid 
(t_1,\dots ,t_{n-1})\in V_{d}(\A^*)\ {\rm and }\ t_1\cdots t_n=1\},
\]
and so the computation of $V_{d}(\A)$ reduces to that of $V_{d}(\A^*)$,  
see \cite{CScv}.

The second construction associates to an affine arrangement, $\A$, 
in $\C^2$, a linearly fibered arrangement, $\wA$, also in $\C^2$, 
called a {\em big arrangement} associated to $\A$. 
The construction depends on the choice of a linear projection 
$\bar{p}:\C^2\to\C$, for which no line of $\A$ coincides with 
$\bar{p}^{-1}(\rm{point})$.  
Let $H_1,\dots, H_n$ be the lines of $\A$, 
let $v_1,\dots,v_s$ be their intersection points, and let 
$\{w_1,\dots,w_r\}=\bar{p}(\{v_1,\dots,v_s\})$. Then 
$\wA=\A\cup \{H_{n+1},\dots, H_{n+r}\}$,  where $H_{n+j}=\bar{p}^{-1}(w_j)$.  
The restriction 
$\bar{p}:\widehat{X}\to\C\setminus\{w_1,\dots ,w_r\}$ is a (linear) fibration, 
with fiber $\C\setminus \{n\ \rm{ points}\}$.  The monodromy generators, 
$\bar{\a}_1,\dots ,\bar{\a}_r$ may be found using a slight modification of the
algorithm from \cite{CSbm} (see also \cite{Coh}).  Deform $\bar{p}$ to a generic
projection
$p:\C^2\to\C$, and let $\a_1,\dots ,\a_s$ be the corresponding braid monodromy
generators.  Then, $\bar{\a}_j=\prod_{p(v_k)=w_j} \a_k$.  The fundamental
group of $\wA$ is the semidirect product
$\widehat{G}=\F_{n}\rtimes_{\bar\a} \F_{r}$, with presentation
\begin{equation}  
\label{eq:linfib}
\widehat{G}=\langle x_1,\dots ,x_n,y_1,\dots ,y_r\mid
x_i^{y_j}=\bar{\a}_j(x_i)\rangle.
\end{equation}
Given an arrangement $\B$ so that $\B=\wA$, the presentation 
(\ref{eq:linfib}) of $\pi_1(X(\wA))=\widehat{G}$ is often simpler to use 
than the presentation (\ref{eq:bmpres}), obtained from the 
general braid monodromy algorithm applied directly to $\B$.  
In particular, if we pick $\{t_1,\dots, t_{n+r}\}=\{x_1,\dots,x_n, 
y_1,\dots, y_{r}\}^{\ab}$ as basis for $H_1(\widehat{G})\cong\Z^{n+r}$, 
the Alexander matrix of $\widehat{G}$ has the block form 
\begin{equation}  
\label{eq:bigal}
A=  
\left( \begin{array}{cccc}
\id -t_{n+1}\cdot \Theta (\bar\a _{1}) & d_{1} &\cdots & 0 \\ 
\vdots & & \ddots \\ 
\id -t_{n+r}\cdot \Theta (\bar\a _{r}) & 0 &\cdots & d_{1}
\end{array}\right) ,
\end{equation}
where $\Theta:P_n\to \GL(n,\Z[t_1^{\pm 1},\dots ,t_n^{\pm 1}])$ 
is the Gassner representation, and 
$d_1=\left( t_{1}-1 \ \cdots \ t_{n}-1 \right)^{\top}$, 
see \cite[\S3.9]{CSai}.

\section{Warm-up Examples}
\label{sec:warmup}

We continue with some relatively simple examples of hyperplane arrangements 
and their characteristic varieties.  These examples, which illustrate
the above discussion, will be useful in understanding subsequent, 
more complicated examples.


\begin{example}
\label{ex:braid}  
{\rm 
Let $\A_3$ be the braid arrangement in $\C^3$,   
with defining polynomial $Q=xyz(x-y)(x-z)(y-z)$. 
The decone $\A_3^*$, obtained by setting $z=1$, 
is depicted in Figure~\ref{fig:dbraid}.  Note that 
$\A_3^*=\wA$, where $\A$ consists of the lines marked $1, 2, 3$.   
Thus, $\A_3$ is fiber-type, with exponents $\{1,2,3\}$, and 
$G^*=\F_3\rtimes_{\bar\a} \F_2$, where $\bar\a_1=A_{12}$, $\bar\a_2=A_{13}$ 
(of course, $G=P_4\cong G^*\times \Z$). 

The resonance and 
characteristic varieties of $\A_3$ were computed in \cite{Fa}, 
\cite{Li3}, \cite{CScv}.  The variety $V_{1}(\A_3)\subset (\C^*)^6$ has 
$4$ local components, corresponding to the triple points $124, 135, 236, 456$, 
and one essential component, corresponding to the neighborly 
partition $(16| 25 | 34)$:
\[
\Pi=\{(s,t,(st)^{-1},(st)^{-1},t,s) \mid s,t\in \C^*\},
\]
see Figure~\ref{fig:mbraid}.  The components of $V_1$ meet only at 
$\b{1}$.  Moreover, $V_{2}=\cdots=V_{6}=\{\b{1}\}$.  
The intersection poset of the characteristic varieties of $\A_3$ 
is depicted in Figure~\ref{fig:cvbraid}.   The poset is ranked 
by dimension (indicated by relative height), and filtered according 
to depth in the characteristic tower  (indicated by color: 
$V_1$ in black, $V_2$ in white).
}
\end{example}

\begin{figure}
\subfigure[Decone $\A_3^*$]{%
\label{fig:dbraid}%
\begin{minipage}[t]{0.21\textwidth}
\setlength{\unitlength}{14pt}
\begin{picture}(4.5,2.5)(-1.1,0)
\multiput(0,1)(0,2){2}{\line(1,0){4}}
\multiput(1,0)(2,0){2}{\line(0,1){4}}
\put(0,4){\line(1,-1){4}}
\put(4.5,0){\makebox(0,0){$1$}}
\put(4.5,1){\makebox(0,0){$2$}}
\put(4.5,3){\makebox(0,0){$3$}}
\put(3,4.5){\makebox(0,0){$4$}}
\put(1,4.5){\makebox(0,0){$5$}}
\end{picture}
\end{minipage}
}
\subfigure[Matroid 
and torus $\Pi$]{%
\label{fig:mbraid}%
\setlength{\unitlength}{0.5cm}
\begin{minipage}[t]{0.37\textwidth}
\begin{picture}(5,3)(-3,-0.5)
\put(3,3){\line(1,-1){3}}
\put(3,3){\line(-1,-1){3}}
\put(1.5,1.5){\line(3,-1){4.5}}
\put(4.5,1.5){\line(-3,-1){4.5}}
\multiput(0,0)(6,0){2}{\circle*{0.4}}
\multiput(1.5,1.5)(3,0){2}{\circle*{0.4}}
\multiput(3,3)(0,-2){2}{\circle*{0.4}}
\multiput(0,-0.6)(6,0){2}{\makebox(0,0){$(st)^{-1}$}}
\multiput(1,1.6)(4,0){2}{\makebox(0,0){$t$}}
\multiput(3,3.5)(0,-3){2}{\makebox(0,0){$s$}}
\end{picture}
\end{minipage}
}
\setlength{\unitlength}{0.8cm}
\subfigure[Characteristic varieties]{
\label{fig:cvbraid}%
\begin{minipage}[t]{0.36\textwidth}
\begin{picture}(3,2.5)(-1.5,-0.4)
\xygraph{!{0;<6mm,0mm>:<0mm,16mm>::}
[]*DL(2){\mathbf{1}}*\cir<3pt>{}
(-[dll]*-{\blacktriangle}
,-[dl]*-{\blacktriangle}
,-[d] *-{\blacktriangle}
,-[dr]*-{\blacktriangle}  
,-[drr]*U(3){\Pi}*-{\blacklozenge})
}
\end{picture}
\end{minipage}
}
\caption{\textsf{The braid arrangement $\A_3$}}
\label{fig:braid}
\end{figure}


%
\begin{figure}
\subfigure[Decone $\NN^*$]{%
\label{fig:dnonfano}%
\begin{minipage}[t]{0.4\textwidth}
\setlength{\unitlength}{14pt}
\begin{picture}(6,5)(-3.7,0)
\multiput(0,1)(0,2){2}{\line(1,0){4}}
\multiput(1,0)(2,0){2}{\line(0,1){4}}
\put(0,4){\line(1,-1){4}}
\put(0,0){\line(1,1){4}}
\put(4.5,0){\makebox(0,0){$1$}}
\put(4.5,1){\makebox(0,0){$2$}}
\put(4.5,3){\makebox(0,0){$3$}}
\put(4.5,4.5){\makebox(0,0){$4$}}
\put(3,4.5){\makebox(0,0){$5$}}
\put(1,4.5){\makebox(0,0){$6$}}
\end{picture}
\end{minipage}
}
\setlength{\unitlength}{0.8cm}
\subfigure[Characteristic varieties of $\NN$]{%
\label{fig:cvnonfano}%
\begin{minipage}[t]{0.45\textwidth}
\begin{picture}(6,2.5)(-1.2,-0.5)
\xygraph{!{0;<6mm,0mm>:<0mm,16mm>::}
[]*DL(2){\mathbf{1}}*\cir<3pt>{}
(-[dllll]*-{\blacktriangle}
,-[dlll]*-{\blacktriangle}
,-[dll]*-{\blacktriangle}
,-[dl]*-{\blacktriangle}
,-[d] *-{\blacktriangle}
,-[dr]*-{\blacktriangle}  
,-[drr]*U(2.8){\Pi_1}*-{\blacklozenge}
,-[drrr]*U(2.8){\Pi_2}*-{\blacklozenge}
,-[drrrr]*U(2.8){\Pi_3}*-{\blacklozenge}
,[rrr]*DL(2){\rho}*\cir<2pt>{}
(-[dl]
,-[d]
,-[dr])
}
\end{picture}
\end{minipage}
}
\caption{\textsf{The non-Fano arrangement $\NN$}}
\label{fig:diamond}
\end{figure}

\begin{example}
\label{ex:diamond}
{\rm 
A realization of the non-Fano plane is the arrangement $\NN$, with  
defining polynomial $Q=xyz(x-y)(x-z)(y-z)(x+y-z)$.  A decone $\NN^*$ 
is depicted in Figure~\ref{fig:dnonfano}.  

The characteristic varieties of $\NN$ were computed in \cite{CScv} 
(see also \cite{LY2}).  
The variety $V_{1}\subset (\C^*)^7$ has $6$ local 
components, corresponding to triple points, and $3$ non-local components, 
$\Pi_1=\Pi(25|36|47)$, $\Pi_2=\Pi(17|26|35)$, $\Pi_3=\Pi(14|23|56)$, 
corresponding to braid sub-arrangements.  The local 
components meet only at $\b{1}$, but the non-local 
components also meet at the point $\rho=(1,-1,-1,1,-1,-1,1)$.  
The variety $V_{2}=\{\b{1},\rho\}$ is a discrete algebraic 
subgroup of $(\C^*)^7$, isomorphic to $\Z_2$.   The characteristic 
intersection poset of $\NN$ is depicted in Figure~\ref{fig:cvnonfano}. 
}
\end{example}


\begin{example}
\label{ex:B3}
{\rm 
Let $\B_3$ be the reflection arrangement of type ${\rm B}_3$,  
with defining polynomial $Q=xyz(x-y)(x-z)(y-z)(x-y-z)(x-y+z)(x+y-z)$.  
A decone is shown in Figure~\ref{fig:B3-a}.  Note that 
$\B_3^*=\wA$, where $\A$ consists of the lines marked $1,\dots, 5$.  Thus, 
$\B_3$ is fiber-type, with exponents $\{1,3,5\}$, and 
$G^*=\F_5\rtimes_{\bar\a} \F_3$, where $\bar\a_1=A_{234}$,
$\bar\a_2=A_{14}^{A_{24}A_{34}}A_{25}$, $\bar\a_3=A_{35}^{A_{23}A_{25}}$. 

A computation with Fox derivatives, using the techniques from 
\S\ref{sec:cvar}, shows that the characteristic variety 
$V_{1}\subset (\C^*)^9$ has $19$ components: 
\begin{itemize}
\item  $7$ local components, corresponding to 
$4$ triple points  and $3$ quadruple points.
\item  $11$ components corresponding to 
braid sub-arrangements.
\item  $1$ essential, $2$-dimensional component, corresponding to the 
neighborly partition $(156 | 248 | 379)$, identified in \cite{Fa}, Example 4.6: 
\[
\Gamma=\{ (t,s,(st)^{-2},s,t,t^2,(st)^{-1},s^2,(st)^{-1}) \mid s, t\in \C^*\}.
\]
\end{itemize}

Three triples of braid components meet $\Gamma$ on $V_{2}$, at the points 
\[
\rho_1=(1,-1,1,-1, 1,1,-1,1,-1),\quad \rho_2=(-1,1,1,1,-1,1,-1,1,-1),
\] 
and $\rho_1\rho_2$.  The variety $V_2$ consists of three 
$3$-dimensional tori (corresponding to  quadruple points), 
together with the discrete subgroup 
$\Z_2^2=\{1,\rho_1,\rho_2,\rho_1\rho_2\}$.
The characteristic 
intersection poset of $\B_3$ is depicted in Figure~\ref{fig:B3-b}. 
}
\end{example}

\begin{figure}
\subfigure[Decone $\B_3^*$]{%
\label{fig:B3-a}%
\begin{minipage}[t]{0.3\textwidth}
\setlength{\unitlength}{0.6cm}
\begin{picture}(5,5)(-0.5,-1)
\multiput(1,0)(1,0){2}{\line(1,1){3}}
\multiput(4,0)(1,0){2}{\line(-1,1){3}}
\multiput(2.5,0)(0.5,0){3}{\line(0,1){3}}
\put(1,1.5){\line(1,0){4}}
\put(4.2,-0.5){\makebox(0,0){$1$}}
\put(5.2,-0.5){\makebox(0,0){$2$}}
\put(5.5,1.5){\makebox(0,0){$3$}}
\put(5.2,3.5){\makebox(0,0){$4$}}
\put(4.2,3.5){\makebox(0,0){$5$}}
\put(3.5,3.5){\makebox(0,0){$6$}}
\put(3,3.5){\makebox(0,0){$7$}}
\put(2.5,3.5){\makebox(0,0){$8$}}
\end{picture}
\end{minipage}
}
\subfigure[Characteristic varieties of $\B_3$]{%
\label{fig:B3-b}%
\begin{minipage}[t]{0.6\textwidth}
\setlength{\unitlength}{0.6cm}
\begin{picture}(6,5)(0,-3.5)
\xygraph{!{0;<4.5mm,0mm>:<0mm,7mm>::}
[]*DL(2){\mathbf{1}}*\cir<3pt>{}
(-[ddddllllllllll]*{\square}
,-[ddddlllllllll]*{\square}
,-[ddddllllllll]*{\square}
,-[dddlllll]*-{\blacktriangle}
,-[dddllll]*-{\blacktriangle}
,-[dddlll]*-{\blacktriangle}
,-[dddll]*-{\blacktriangle} 
,-[dddl]*-{\blacklozenge}
,-[ddd]*-{\blacklozenge}
,-[dddr]*-{\blacklozenge}
,-[dddrr]*-{\blacklozenge}
,-[dddrrr]*-{\blacklozenge}
,-[dddrrrr]*-{\blacklozenge}
,-[dddrrrrr]*-{\blacklozenge}
,-[dddrrrrrr]*-{\blacklozenge}
,-[dddrrrrrrr]*-{\blacklozenge}
,-[dddrrrrrrrr]*-{\blacklozenge}
,-[dddrrrrrrrrr]*-{\blacklozenge}
,-[dddrrrrrrrrrr]*U(2.5){\Gamma}*-{\blacksquare}
,[rrrr]*DL(2){\rho_1}*\cir<2pt>{}
(-[dddl]
,-[ddd]
,-[dddrr]
,-[dddrrrrrr])
,[rrrrrr]*DL(2){\rho_2}*\cir<2pt>{}  
(-[dddl]
,-[ddd]
,-[dddr]
,-[dddrrrr])
,[rrrrrrrr]*DL(1.5){\rho_1\rho_2}*\cir<2pt>{}  
(-[dddll]
,-[ddd]
,-[dddr]
,-[dddrr])
)
}
\end{picture}
\end{minipage}
}
\caption{\textsf{The ${\rm B}_3$ reflection arrangement}}
\label{fig:B3arr}
\end{figure}

\section{Positive-dimensional translated tori}
\label{sec:posdim}
We now come to our basic example of a complex hyperplane arrangement 
whose top characteristic variety contains a positive-dimensional 
translated component.   


\begin{example}
\label{ex:deletedb3}

\begin{figure}[t]
\subfigure[Decone $\DD^*$]{%
\label{fig:delB3-a}%
\begin{minipage}[t]{0.33\textwidth}
\setlength{\unitlength}{0.75cm}
\begin{picture}(5,5)(-0.3,-0.6)
\multiput(1,0)(1,0){2}{\line(1,1){3}}
\multiput(4,0)(1,0){2}{\line(-1,1){3}}
\multiput(2.5,0)(0.5,0){3}{\line(0,1){3}}
\put(4.2,-0.5){\makebox(0,0){$1$}}
\put(5.2,-0.5){\makebox(0,0){$2$}}
\put(5.2,3.5){\makebox(0,0){$3$}}
\put(4.2,3.5){\makebox(0,0){$4$}}
\put(3.5,3.5){\makebox(0,0){$5$}}
\put(3,3.5){\makebox(0,0){$6$}}
\put(2.5,3.5){\makebox(0,0){$7$}}
\end{picture}
\end{minipage}
}
\setlength{\unitlength}{0.65cm}
\subfigure[Matroid of $\DD$ and parametrization of $C$]{%
\label{fig:delB3-b}%
\begin{minipage}[t]{0.6\textwidth}
\begin{picture}(6,5)(-4,-0.7)
\put(0,2){\line(1,0){6}}
\put(0,0){\line(0,1){4}}
\put(0,0){\line(3,1){6}}
\put(0,4){\line(3,-1){6}}
\put(0,0){\line(1,1){3}}
\put(0,4){\line(1,-1){3}}
\put(3,1){\line(0,1){2}}
\multiput(0,0)(0,2){3}{\circle*{0.4}}
\multiput(3,1)(0,1){3}{\circle*{0.4}}
\put(6,2){\circle*{0.4}}
\put(2,2){\circle*{0.4}}
\put(-0.5,-0.7){\makebox(0,0){$-t^{-1}$}}
\put(-0.6,2){\makebox(0,0){$t^2$}}
\put(-0.5,4.5){\makebox(0,0){$-t^{-1}$}}
\put(1.9,1.4){\makebox(0,0){$-1$}}
\put(3,0.4){\makebox(0,0){$t$}}
\put(3.6,2.4){\makebox(0,0){$t^{-2}$}}
\put(3,3.6){\makebox(0,0){$t$}}
\put(6.8,2){\makebox(0,0){$-1$}}
\end{picture}
\end{minipage}
}
\caption{\textsf{The deleted ${\rm B}_3$ arrangement $\DD$}}
\label{fig:deconeA}
\end{figure}

{\rm
Let $\DD$ be the arrangement obtained from the ${\rm B}_3$ 
reflection arrangement by deleting the plane $x+y-z=0$.  A defining 
polynomial for $\DD$ is $Q=xyz(x-y)(x-z)(y-z)(x-y-z)(x-y+z)$.
The decone $\DD^*$, obtained by setting $z=1$, 
is depicted in Figure~\ref{fig:delB3-a}.  
Note that $\DD^*=\wA$, where $\A$ consists of the lines marked 
$1,\dots, 4$.  Thus,  
$\DD$ is fiber-type, with exponents $\{1,3,4\}$, and 
$G^*=\F_4\rtimes_{\bar\a} \F_3$, where $\bar\a_1=A_{23}$,
$\bar\a_2=A_{13}^{A_{23}}A_{24}$,  $\bar\a_3=A_{14}^{A_{24}}$. 

The Alexander matrix of $G^*$, given by (\ref{eq:bigal}), 
is row-equivalent to 
{\scriptsize{
\[
\mbox{{\normalsize $A=$}}
\left(
\begin{array}{ccccccc}
1-t_{5}&0&0&0&t_{1}-1&0&0\\
0&t_{5}(t_{3}-1)&1-t_{2}t_{5}&0&t_{3}-1&0&0\\
0&1-t_{5}&t_{2}(1-t_{5})&0&t_{2}t_{3}-1&0&0\\
0&0&0&1-t_{5}&t_{4}-1&0&0\\
t_{6}(t_{3}-1)&(t_{3}-1)(t_{1}t_{6}-1)&t_{2}(1-t_{1}t_{6})&0&0&t_{3}-1&0\\
1-t_{6}&t_{1}(t_{3}-1)(t_{6}-1)&t_{1}t_{2}(1-t_{6})&0&0&t_{1}t_{3}-1&0\\
0&t_{6}(t_{4}-1)&0&1-t_{2}t_{6}&0&t_{4}-1&0\\
0&1-t_{6}&0&t_{2}(1-t_{6})&0&t_{2}t_{4}-1&0\\
t_{7}(t_{4}-1)&(t_{4}-1)(t_{1}t_{7}-1)&0&t_{2}(1-t_{1}t_{7})&0&0&t_{4}-1\\
1-t_{7}&t_{1}(1-t_{7})&0&t_{1}t_{2}(1-t_{7})&0&0&t_{1}t_{2}t_{4}-1\\
0&1-t_{7}&0&0&0&0&t_{2}-1\\
0&0&1-t_{7}&(t_{3}-1)(1-t_{7})&0&0&t_{4}(t_{3}-1)
\end{array}
\right).  
\]
}}
Now recall that $V_{1}(\DD)=\{\mathbf{t}\in (\C^*)^8 \mid (t_1,\dots,t_7)\in
V_1(\DD^*)\ {\rm and }\ t_1\cdots t_8=1\}$, where 
$V_1(\DD^*)$ is the sub-variety of $(\C^*)^7$ defined by 
the ideal of $6\times 6$ minors of the matrix $A$.   Computing the primary 
decomposition of that ideal reveals that the variety $V_{1}(\DD)$ has
$13$ components:
\begin{itemize}
\item  $7$ local components, corresponding to $6$ triple points and 
one quadruple point.  
\item  $5$ non-local components passing through $\b{1}$, 
corresponding to braid sub-arrangements: 
$\Pi_1=\Pi(15 | 26 | 38)$, 
$\Pi_2=\Pi(28 | 36 | 45)$, 
$\Pi_3=\Pi(14 | 23 | 68)$, 
$\Pi_4=\Pi(16 | 27 | 48)$, 
$\Pi_5=\Pi(18 | 37 | 46)$. 
\item  $1$ essential component, which does not pass through $\mathbf{1}$.  
This component is $1$-dimensional, and is parametrized by
\[
C=\{ (t,-t^{-1},-t^{-1},t,t^2,-1,t^{-2},-1) \mid t\in \C^*\}.
\]
\end{itemize}
(Note that the translated torus $C$ is one of the two connected 
components of $\Gamma\cap \{t_3=1\}$.) The braid components of 
$V_1(\DD)$ meet $C$ at the points 
\[
\begin{array}{ll}
\rho_1&=\Pi_1\cap\Pi_2\cap\Pi_3\cap C=(1,-1,-1,1,1,-1,1,-1),\\
\rho_2&=\Pi_3\cap\Pi_4\cap\Pi_5\cap C=(-1,1,1,-1,1,-1,1,-1),
\end{array}
\]
both of which belong to $V_{2}(\DD)$.  
The characteristic 
intersection poset of $\DD$ is depicted in Figure~\ref{fig:cva}. 
}
\end{example}

\begin{figure}
\setlength{\unitlength}{12pt}
\begin{picture}(12,8)(-5,-5)
\xygraph{!{0;<7mm,0mm>:<0mm,7mm>::}
[]*DL(2){\mathbf{1}}*\cir<3pt>{}
(-[ddddllllll]*{\square}
,-[dddllll]*-{\blacktriangle}
,-[dddlll]*-{\blacktriangle}
,-[dddll]*-{\blacktriangle}
,-[dddl]*-{\blacktriangle}
,-[ddd] *-{\blacktriangle}
,-[dddr]*-{\blacktriangle}  
,-[dddrr]*U(2.5){\Pi_1}*-{\blacklozenge}
,-[dddrrr]*U(2.5){\Pi_2}*-{\blacklozenge}
,-[dddrrrr]*U(2.5){\Pi_3}*-{\blacklozenge}
,-[dddrrrrr]*U(2.5){\Pi_4}*-{\blacklozenge}
,-[dddrrrrrr]*U(2.5){\Pi_5}*-{\blacklozenge}
,[drrrr]*U(2.5){C}*{\bigstar}
,[rrr]*DL(2){\rho_1}*\cir<2pt>{}
(-[dddl]
,-[ddd]
,-[dddr]
,-[dr])
,[rrrrr]*DL(2){\rho_2}*\cir<2pt>{}  
(-[dddl]
,-[ddd]
,-[dddr]
,-[dl])
)
}
\end{picture}
\caption{\textsf{The characteristic varieties of $\DD$}}
\label{fig:cva}
\end{figure}

As noted in the Introduction, this example answers Problem 5.1
in \cite{LY2}.  It also answers Problems 5.2 and 5.3 in \cite{LY2}.  
Indeed, let $\lambda=(\frac{1}{4},\frac{1}{4},\frac{1}{4},\frac{1}{4},
\frac{1}{2},-\frac{1}{2},-\frac{1}{2},-\frac{1}{2})$.  
Clearly, $\lambda$, and all its integral translates, do 
not belong to $\RR_1(\DD)$, because all components of 
$\RR_1(\DD)$ are non-essential.  Hence, $H^1(H^*(X,\C),\lambda+N)=0$,  
for all $N\in\Z^8$.  On the other hand, 
$\b{t}:=\exp(\lambda)= (\ii,\ii,\ii,\ii,-1,-1,-1,-1)$ belongs 
to $C$, and thus $\dim H^1(X,\C_{\b{t}})=1$. 

Finally, this example also answers in the negative 
Conjecture 4.4 from \cite{LY1}, at least in its strong form.  Indeed, 
there are infinitely many $\b{t}=\exp (\lambda)\in C$ for which 
\[
1=\dim H^1(X,\C_{\b{t}})\ne \sup_{N\in \Z^8} \dim H^1(H^*(X,\C),\lambda+N)=0.
\]

\section{Further examples}
\label{sec:further}

In this section, we give a few more examples that illustrate the nature of 
translated components in the characteristic varieties of arrangements. 


\begin{example}  
\label{ex:grunbaum}
{\rm
Let $\A={\rm A}_2(10)$ be the simplicial arrangement from 
the list in Gr\"{u}nbaum~\cite{Gr}. A defining polynomial for $\A$ is 
$Q=xyz(y-x)(y+x)(2y-z)(y-x-z)(y-x+z)(y+x+z)(y+x-z)$.  
Figure~\ref{fig:grunbaum} shows a decone $\A^*$, 
together with the characteristic intersection poset of $\A$. 
The variety $V_1(\A)$ has $33$ components:
\begin{itemize}
\item  $10$ local components, corresponding to 
$7$ triple and $3$ quadruple points.
\item  $17$ non-local components corresponding
to braid sub-arrangements.
\item  $1$ non-local component, $\Gamma=\Gamma(1 5 10 | 2 4 8 | 3 7 9)$,  
corresponding to a $\B_3$ sub-arrangement.
\item  $3$ components that do not pass through the origin,  
$C_1=C(1  3  4  5  6  8  9  10)$, 
$C_2=C(2  3  4  5  6  8  9  10)$, 
$C_3=C(3  4  5  6  7  8  9  10)$, 
corresponding to $\DD$ sub-arrangements.
\item  $2$ isolated points of order $6$,  
$\zeta=(\eta^2,\eta^2,\eta,\eta^2,\eta^2,-1,\eta^2,\eta,\eta^2,\eta)$ 
and $\zeta^{-1}$, where $\eta=e^{\pi\ii/3}$. 
\end{itemize}

Note that all positive-dimensional components of $V_1(\A)$ 
are non-essential, whereas the two $0$-dimensional components 
are essential.  The non-local components meet at $7$ isolated 
points of order~$2$, belonging to $V_2(\A)$: 
\[
\hspace*{-0.2in}
\begin{array}{lll}
\rho_1&=(1, 1, -1, -1, 1, 1, 1, -1, -1, 1), \
\rho_2&=(1, 1, -1, 1, -1, 1, 1, 1, -1, -1),
\\
\rho_3&=(-1, -1, 1, -1, -1, 1, 1, 1, 1, 1),\ 
\rho_4&=(-1, 1, 1, 1, -1, 1, -1, 1, -1, 1),
\end{array}
\]
$\rho_1\rho_2$, $\rho_3\rho_4$, and $\zeta^3$. 
}

\begin{figure}
\subfigure{
\label{fig:grundec}%
\begin{minipage}[t]{0.23\textwidth}
\setlength{\unitlength}{0.62cm}
\begin{picture}(4,6)(0.5,-2)
\put(1,0.5){\line(1,1){2.75}}
\put(1.5,0){\line(1,1){3}}
\put(2.25,-0.25){\line(1,1){2.75}}
\put(3.75,-0.25){\line(-1,1){2.75}}
\put(4.5,0){\line(-1,1){3}}
\put(5,0.5){\line(-1,1){2.75}}
\put(3,-0.5){\line(0,1){4}}
\put(1,1.5){\line(1,0){4}}
\put(1,2){\line(1,0){4}}
\put(3,-1){\makebox(0,0){$1$}}
\put(3.9,-0.8){\makebox(0,0){$2$}}
\put(4.8,-0.5){\makebox(0,0){$3$}}
\put(5.5,0.2){\makebox(0,0){$4$}}
\put(5.5,1.4){\makebox(0,0){$5$}}
\put(5.5,2.1){\makebox(0,0){$6$}}
\put(5.5,2.9){\makebox(0,0){$7$}}
\put(4.8,3.4){\makebox(0,0){$8$}}
\put(4,3.7){\makebox(0,0){$9$}}
\end{picture}
\end{minipage}
}
\subfigure{
\label{fig:gruncv}%
\begin{minipage}[t]{0.67\textwidth}
\setlength{\unitlength}{0.45cm}
\begin{picture}(5,6)(0.5,-6)
\xygraph{!{0;<3.7mm,0mm>:<0mm,7mm>::}
[]*DL(2){\mathbf{1}}*\cir<3pt>{}
(-[dddddllllllll]*{\square}
,-[dddddlllllll]*{\square}
,-[dddddllllll]*{\square}
,-[ddddllll]*-{\blacktriangle}
,-[ddddlll]*-{\blacktriangle} 
,-[ddddll]*-{\blacktriangle}
,-[ddddl]*-{\blacktriangle} 
,-[dddd]*-{\blacktriangle} 
,-[ddddr]*-{\blacktriangle}
,-[ddddrr]*-{\blacktriangle} 
,-[ddddrrr]*-{\blacklozenge}
,-[ddddrrrr]*-{\blacklozenge}
,-[ddddrrrrr]*-{\blacklozenge}
,-[ddddrrrrrr]*-{\blacklozenge}
,-[ddddrrrrrrr]*-{\blacklozenge}
,-[ddddrrrrrrrr]*-{\blacklozenge}
,-[ddddrrrrrrrrr]*-{\blacklozenge}
,-[ddddrrrrrrrrrr]*-{\blacklozenge}
,-[ddddrrrrrrrrrrr]*-{\blacklozenge}
,-[ddddrrrrrrrrrrrr]*-{\blacklozenge}
,-[ddddrrrrrrrrrrrrr]*-{\blacklozenge}
,-[ddddrrrrrrrrrrrrrr]*-{\blacklozenge}
,-[ddddrrrrrrrrrrrrrrr]*-{\blacklozenge}
,-[ddddrrrrrrrrrrrrrrrr]*-{\blacklozenge}
,-[ddddrrrrrrrrrrrrrrrrr]*-{\blacklozenge}
,-[ddddrrrrrrrrrrrrrrrrrr]*-{\blacklozenge}
,-[ddddrrrrrrrrrrrrrrrrrrr]*-{\blacklozenge}
,-[ddddrrrrrrrrrrrrrrrrrrrr]*U(2.5){\Gamma}*-{\blacksquare}
,[drrrrrr]*D(2.2){^{C_1}}*{\bigstar}
,[drrrrrrrr]*D(2.4){^{C_2}}*{\bigstar}
,[drrrrrrrrrr]*D(2.2){^{C_3}}*{\bigstar}
,[rrrrr]*D(2){\rho_1}*\cir<2pt>{} 
(-[dr]
,-[ddddl]
,-[dddd]
,-[ddddr])
,[rrrrrrr]*D(2){\zeta^3}*\cir<2pt>{}  
(:@/_/@{-} [drrr]
,-[dl]
,-[dr]
,-[ddddl]
,-[dddd]
,-[ddddr])
,[rrrrrrrrr]*D(2){\rho_2}*\cir<2pt>{} 
(-[dl]
,-[ddddl]
,-[dddd]
,-[ddddr])
,[rrrrrrrrrrr]*D(2){\rho_1\rho_2}*\cir<2pt>{} 
(:@/^/@{-} [ddddllll]
,-[dl]
,-[dddd]
,-[ddddr])
,[rrrrrrrrrrrrrr]*D(2){\rho_3}*\cir<2pt>{} 
(-[ddddl]
,-[dddd]
,-[ddddr]
,-[ddddrrrrrr])
,[rrrrrrrrrrrrrrrr]*D(2){\rho_4}*\cir<2pt>{} 
(-[ddddl]
,-[dddd]
,-[ddddr]
,-[ddddrrrr])
,[rrrrrrrrrrrrrrrrrr]*D(2){\rho_3\rho_4}*\cir<2pt>{}  
(-[ddddlll]
,-[dddd]
,-[ddddr]
,-[ddddrr])
,[ll]*D(2){\zeta^{-1}}*{\bullet}
,[llll]*D(2){\zeta}*{\bullet}
)
}
\end{picture}
\end{minipage}
}
\caption{\textsf{The Gr\"{u}nbaum arrangement 
${\rm A}_2(10)$}}
\label{fig:grunbaum}
\end{figure}

\end{example} 


\begin{example}
\label{ex:falk}
{\rm
Consider the arrangements $\FF_1$ and $\FF_2$, 
with defining polynomials
$Q_1=(x-y-z)Q$ and $Q_2=(x-y-2z)Q$, where 
$Q=xyz(x-y)(y-z)(x-z)(x-2z)(x-3z)$.  
This pair of arrangements was introduced by Falk in \cite{Fa1}.   
Their decones and characteristic varieties 
are depicted in Figure~\ref{fig:falk}. 

\begin{figure}
\subfigure{
\label{fig:decf1}%
\begin{minipage}[t]{0.18\textwidth}
\setlength{\unitlength}{0.4cm}
\begin{picture}(3,6)(0.5,0)
\multiput(2,2)(0,-1){2}{\line(1,1){4}}
\multiput(2.5,0)(1,0){4}{\line(0,1){6}}
\multiput(2,2.5)(0,1){2}{\line(1,0){4}}
\end{picture}
\end{minipage}
}
\subfigure{
\label{fig:cvf1}%
\begin{minipage}[t]{0.29\textwidth}
\setlength{\unitlength}{10pt}
\begin{picture}(6,5)(0.2,-5)
\xygraph{!{0;<3mm,0mm>:<0mm,6mm>::}
[]*DL(2){\mathbf{1}}*\cir<3pt>{}
(-[ddddlllllll]*{\boxplus}
,-[dddllll]*-{\blacktriangle}
,-[dddlll]*-{\blacktriangle}
,-[dddll]*-{\blacktriangle}
,-[dddl]*-{\blacktriangle}
,-[ddd] *-{\blacktriangle}
,-[dddr]*-{\blacktriangle}  
,-[dddrr]*-{\blacklozenge}
,-[dddrrr]*-{\blacklozenge}
,-[dddrrrr]*-{\blacklozenge}
,-[dddrrrrr]*-{\blacklozenge}
,-[dddrrrrrr]*-{\blacklozenge}
,[drrrr]*{\bigstar}
,[rrr]*\cir<2pt>{}
(-[dddl]
,-[ddd]
,-[dddr]
,-[dr])
,[rrrrr]*\cir<2pt>{}  
(-[dddl]
,-[ddd]
,-[dddr]
,-[dl])
)
}
\end{picture}
\end{minipage}
}
\subfigure{
\label{fig:decf2}%
\begin{minipage}[t]{0.17\textwidth}
\setlength{\unitlength}{0.4cm}
\begin{picture}(3,6)(0,-1)
\multiput(2,1)(0,-2){2}{\line(1,1){4}}
\multiput(2.5,-1)(1,0){4}{\line(0,1){6}}
\multiput(2,1.5)(0,1){2}{\line(1,0){4}}
\end{picture}
\end{minipage}
}
\subfigure{
\label{fig:cvf2}%
\begin{minipage}[t]{0.29\textwidth}
\setlength{\unitlength}{10pt}
\begin{picture}(6,5)(-0.5,-5)
\xygraph{!{0;<3mm,0mm>:<0mm,6mm>::}
[]*DL(2){\mathbf{1}}*\cir<3pt>{}
(-[ddddlllllll]*{\boxplus}
,-[dddllll]*-{\blacktriangle}
,-[dddlll]*-{\blacktriangle}
,-[dddll]*-{\blacktriangle}
,-[dddl]*-{\blacktriangle}
,-[ddd] *-{\blacktriangle}
,-[dddr]*-{\blacktriangle}  
,-[dddrr]*-{\blacklozenge}
,-[dddrrr]*-{\blacklozenge}
,-[dddrrrr]*-{\blacklozenge}
,-[dddrrrrr]*-{\blacklozenge}
)
}
\end{picture}
\end{minipage}
}
\caption{\textsf{The Falk fiber-type arrangements $\FF_1$ and $\FF_2$}}
\label{fig:falk}
\end{figure}

Both arrangements are fiber-type, with exponents $\{1, 4, 4\}$.  
Thus, by the LCS formula of Falk and Randell (see \cite{OT}), 
the ranks, $\phi_k(G):=\rank \gr_k(G)$, of the lower central 
series  quotients of the two groups  
are the same.  As noted in \cite{CSpn}, though, 
the ranks, $\theta_k(G):=\rank \gr_k(G/G'')$, of 
the Chen groups are different:
$\theta_k(G(\FF_1)) = \frac{1}{2}(k-1)(k^2+3k+24)$ and 
$\theta_k(G(\FF_2)) = \frac{1}{2}(k-1)(k^2+3k+22)$, for $k\ge 4$.  
Moreover, as noted in \cite{Fa}, the resonance varieties of the two 
arrangements are not isomorphic, even as abstract varieties:  
$\RR_{1}(\FF_1)$ has $12$ components, whereas  $\RR_{1}(\FF_2)$ 
has $11$ components.   
An even more pronounced difference
shows up in the  characteristic varieties:  $V_{1}(\FF_1)$ has a 
$13^{\rm{th}}$ component (corresponding to a sub-arrangement 
isomorphic to $\DD$), which does not pass through the origin.
}

\end{example}


\begin{example}
\label{ex:ziegler}
{\rm
Consider the arrangements $\ZZ_1$ and $\ZZ_2$, 
with defining polynomials 
$Q_1=(x-y-2z)Q$ and $Q_2=(x-y-3z)Q$, where 
$Q=xyz(x-y)(y-z)(x-z)(x-2z)(x-3z)(x-4z)(x-5z)(x-y-z)(x-y-4z)$. 
This pair of arrangements was introduced by Ziegler~\cite{Zi}.  
Their decones and characteristic varieties are depicted in 
Figures~\ref{fig:ziegler1} and \ref{fig:ziegler2}.   

\begin{figure}
\subfigure{
\label{fig:decz1}%
\begin{minipage}[t]{0.22\textwidth}
\setlength{\unitlength}{0.38cm}
\begin{picture}(4,9)(1,-2.5)
\multiput(2,1)(0,-1){3}{\line(1,1){6}}
\put(2,-3){\line(1,1){6}}
\multiput(2.5,-3)(1,0){6}{\line(0,1){10}}
\multiput(2,1.5)(0,1){2}{\line(1,0){6}}
\end{picture}
\end{minipage}
}
\subfigure{
\label{fig:cvz1}%
\begin{minipage}[t]{0.62\textwidth}
\setlength{\unitlength}{0.7cm}
\begin{picture}(6,4)(0.8,-4)
\xygraph{!{0;<3.05mm,0mm>:<0mm,6.2mm>::}
[]*DL(2){\mathbf{1}}*\cir<3pt>{}
(-[ddddddlllllllllllllllllll]*!{\boxtimes}
,-[dddddllllllllllllll]*!{\boxplus}
,-[ddddllllllllll]*-{\blacktriangle}
,-[ddddlllllllll]*-{\blacktriangle}
,-[ddddllllllll]*-{\blacktriangle}
,-[ddddlllllll]*-{\blacktriangle} 
,-[ddddllllll]*-{\blacktriangle}
,-[ddddlllll]*-{\blacktriangle} 
,-[ddddllll]*-{\blacktriangle} 
,-[ddddlll]*-{\blacktriangle}
,-[ddddll]*-{\blacktriangle} 
,-[ddddl]*-{\blacklozenge}
,-[dddd]*-{\blacklozenge}
,-[ddddr]*-{\blacklozenge}
,-[ddddrr]*-{\blacklozenge}
,-[ddddrrr]*-{\blacklozenge}
,-[ddddrrrr]*-{\blacklozenge}
,-[ddddrrrrr]*-{\blacklozenge}
,-[ddddrrrrrr]*-{\blacklozenge}
,-[ddddrrrrrrr]*-{\blacklozenge}
,-[ddddrrrrrrrr]*-{\blacklozenge}
,-[ddddrrrrrrrrr]*-{\blacklozenge}
,-[ddddrrrrrrrrrr]*-{\blacklozenge}
,-[ddddrrrrrrrrrrr]*-{\blacklozenge}
,-[ddddrrrrrrrrrrrr]*-{\blacklozenge}
,-[ddddrrrrrrrrrrrrr]*-{\blacklozenge}
,-[ddddrrrrrrrrrrrrrr]*-{\blacklozenge}
,-[ddddrrrrrrrrrrrrrrr]*-{\blacklozenge}
,-[ddddrrrrrrrrrrrrrrrr]*-{\blacklozenge}
,[drrrrrrrrrr]*{\bigstar}
,[drrrrrrrrrrrr]*{\bigstar}
,[drrrrrrrrrrrrrr]*{\bigstar}
,[rrrrrrrrr]*\cir<2pt>{}
(-[dr]
,-[ddddl]
,-[dddd]
,-[ddddr])
,[rrrrrrrrrrr]*\cir<2pt>{}  
(-[dl]
,-[dr]
,-[ddddl]
,-[dddd]
,-[ddddr])
,[rrrrrrrrrrrrr]*\cir<2pt>{}  
(-[dl]
,-[dr]
,-[ddddl]
,-[dddd]
,-[ddddr])
,[rrrrrrrrrrrrrrr]*\cir<2pt>{}  
(-[dl]
,-[ddddl]
,-[dddd]
,-[ddddr])
)
}
\end{picture}
\end{minipage}
}
\caption{\textsf{The Ziegler fiber-type arrangement $\ZZ_1$}}
\label{fig:ziegler1}
\end{figure}

\begin{figure}
\subfigure{
\label{fig:decz2}%
\begin{minipage}[t]{0.22\textwidth}
\setlength{\unitlength}{0.38cm}
\begin{picture}(4,9)(1,-2.5)
\multiput(2,1)(0,-1){2}{\line(1,1){6}}
\multiput(2,-2)(0,-1){2}{\line(1,1){6}}
\multiput(2.5,-3)(1,0){6}{\line(0,1){10}}
\multiput(2,1.5)(0,1){2}{\line(1,0){6}}
\end{picture}
\end{minipage}
}
\subfigure{
\label{fig:cvz2}%
\begin{minipage}[t]{0.62\textwidth}
\setlength{\unitlength}{0.7cm}
\begin{picture}(6,4)(0.8,-4)
\xygraph{!{0;<3.05mm,0mm>:<0mm,6.2mm>::}
[]*DL(2){\mathbf{1}}*\cir<3pt>{}
(-[ddddddlllllllllllllllllll]*!{\boxtimes}
,-[dddddllllllllllllll]*!{\boxplus}
,-[ddddllllllllll]*-{\blacktriangle}
,-[ddddlllllllll]*-{\blacktriangle}
,-[ddddllllllll]*-{\blacktriangle}
,-[ddddlllllll]*-{\blacktriangle} 
,-[ddddllllll]*-{\blacktriangle}
,-[ddddlllll]*-{\blacktriangle} 
,-[ddddllll]*-{\blacktriangle} 
,-[ddddlll]*-{\blacktriangle}
,-[ddddll]*-{\blacktriangle} 
,-[ddddl]*-{\blacklozenge}
,-[dddd]*-{\blacklozenge}
,-[ddddr]*-{\blacklozenge}
,-[ddddrr]*-{\blacklozenge}
,-[ddddrrr]*-{\blacklozenge}
,-[ddddrrrr]*-{\blacklozenge}
,-[ddddrrrrr]*-{\blacklozenge}
,-[ddddrrrrrr]*-{\blacklozenge}
,-[ddddrrrrrrr]*-{\blacklozenge}
,-[ddddrrrrrrrr]*-{\blacklozenge}
,-[ddddrrrrrrrrr]*-{\blacklozenge}
,-[ddddrrrrrrrrrr]*-{\blacklozenge}
,-[ddddrrrrrrrrrrr]*-{\blacklozenge}
,-[ddddrrrrrrrrrrrr]*-{\blacklozenge}
,-[ddddrrrrrrrrrrrrr]*-{\blacklozenge}
,-[ddddrrrrrrrrrrrrrr]*-{\blacklozenge}
,-[ddddrrrrrrrrrrrrrrr]*-{\blacklozenge}
,-[ddddrrrrrrrrrrrrrrrr]*-{\blacklozenge}
,[drrrrrrrrr]*{\bigstar}
,[drrrrrrrrrrrrrr]*{\bigstar}
,[rrrrrrrr]*\cir<2pt>{}
(-[dr]
,-[ddddl]
,-[dddd]
,-[ddddr])
,[rrrrrrrrrr]*\cir<2pt>{}  
(-[dl]
,-[ddddl]
,-[dddd]
,-[ddddr])
,[rrrrrrrrrrrrr]*\cir<2pt>{}  
(-[dr]
,-[ddddl]
,-[dddd]
,-[ddddr])
,[rrrrrrrrrrrrrrr]*\cir<2pt>{}  
(-[dl]
,-[ddddl]
,-[dddd]
,-[ddddr])
)
}
\end{picture}
\end{minipage}
}
\caption{\textsf{The Ziegler fiber-type arrangement $\ZZ_2$}}
\label{fig:ziegler2}
\end{figure}

Both arrangements are fiber-type, with exponents $\{1, 6, 6\}$; 
thus, $\phi_k(G(\ZZ_1))=\phi_k(G(\ZZ_2))$.  
Even more, the ranks of the Chen groups 
are the same: $\theta_1=13$, $\theta_2=30$, $\theta_3=140$, 
and $\theta_k=\frac{1}{24}(k-1)(k^4+10k^3+47k^2+86k+696)$, for $k\ge 4$.   
Moreover, $\RR_1(\ZZ_1)\cong \RR_1(\ZZ_2)$ (as varieties), 
although one may show, by a rather long calculation of the 
respective polymatroids, that there is no linear isomorphism 
$\C^{13}\to \C^{13}$ taking $\RR_1(\ZZ_1)$ to $\RR_1(\ZZ_2)$. 

On the other hand, the two groups {\em can} be distinguished numerically 
by their characteristic varieties: $V_{1}(\ZZ_1)$ has $32$ components, 
whereas $V_{1}(\ZZ_2)$ has $31$ components.   
Both varieties have $11$ local components 
(corresponding to $9$ triple points, $1$ quintuple point, 
and $1$ septuple point), and $18$ components corresponding to braid 
sub-arrangements.  In addition, both varieties have components which 
do not pass through $\b{1}$, corresponding to $\DD$ sub-arrangements:   
$V_{1}(\ZZ_1)$ has $3$ such components, $V_{1}(\ZZ_2)$ has only $2$. 
}

\end{example}


\section{Concluding remarks}
\label{sec:conclude}

We conclude with a few questions raised by the above examples.  

Let $\A$ be an arrangement of $n$ complex hyperplanes, and let 
$V_d(\A)\subset (\C^*)^n$ ($1\le d\le n$) be its characteristic 
varieties. 
 
\begin{question}  
\label{question1}
Are the translated components of $V_d(\A)$ 
combinatorially determined?
\end{question}

This problem was posed in \cite{CScv} and \cite{LY2}, before 
the existence of translated components in $V_1(\A)$ was known.  
Recall that $\check{V}_d(\A)$---the union of the 
components of $V_d(\A)$ passing through the identity of the torus 
$(\C^*)^n$---{\em is} combinatorially determined.  
If $\check{V}_d(\A)\ne V_d(\A)$ (as in the examples from 
\S\S\ref{sec:posdim}--\ref{sec:further}), the question 
is whether $V_d(\A)\setminus \check{V}_d(\A)$ is also 
determined by the intersection lattice of $\A$. 

\begin{question} 
\label{question2}
What are the possible dimensions of the 
translated components of $V_d(\A)$?
\end{question}

The (positive-dimensional) components passing through the origin 
must have dimension at least $2$, and all dimensions between $2$ 
and $n-1$ can be realized. On the other hand, at least in the 
examples we gave here, the components not passing through $\b{1}$ 
have dimension either $0$ or $1$.  The question is whether 
$V_d(\A)\setminus \check{V}_d(\A)$ can have higher-dimensional 
components. 

\begin{question} 
\label{question3}
What are the possible orders of translation of the 
components of $V_d(\A)$?
\end{question}

In our examples, the components not passing through the origin 
are translated by characters of order $2$ or $6$.  The question 
is whether other orders of translation can occur.  Furthermore, 
one may ask (as a weak form of Question~\ref{question1}) whether 
the orders of translation are combinatorially  determined.  We 
know of a combinatorial upper bound on the lowest common multiple 
of these orders, but do not know when this bound is attained. 


\ack 

I wish to thank D.~Cohen and D.~Matei for useful discussions. 
The computations for this work were done primarily with {\it Macaulay~2} 
(\cite{GS}) and {\it Mathematica~4.0}.  Additional computations were done 
with {\it GAP~4.1} (\cite{gap}).


\end{document}